\documentclass[10pt,leqno]{amsart}
\usepackage{graphicx}
\usepackage{indentfirst,csquotes}

\usepackage{lipsum}
\usepackage{amsfonts}
\usepackage{amsmath}
\usepackage{amssymb}
\usepackage{graphicx}
\usepackage{epstopdf}
\usepackage{caption}
\usepackage{subcaption}
\usepackage{algorithm}
\usepackage[noend]{algpseudocode}
\usepackage{makecell}
\usepackage{multirow}
\usepackage{booktabs}
\usepackage{hyperref}
\usepackage{cleveref}

\newcommand{\R}{\ensuremath{\mathbb{R}}}

\newcommand{\trans}{\ensuremath{\mkern-1.5mu\mathsf{T}}}

\newcommand{\K}{\ensuremath{\mathbf{K}}}
\newcommand{\C}{\ensuremath{\mathbf{C}}}
\newcommand{\M}{\ensuremath{\mathbf{M}}}
\newcommand{\F}{\ensuremath{\mathbf{f}}}
\newcommand{\G}{\ensuremath{\mathbf{g}}}
\newcommand{\V}{\ensuremath{\mathbf{V}}}
\newcommand{\p}{\ensuremath{\mathbf{p}}}
\newcommand{\x}{\ensuremath{\mathbf{x}}}

\newcommand{\xr}{\ensuremath{\mathbf{x}_r}}

\newcommand{\Krk}{\ensuremath{\mathbf{K}_{r,k}}}
\newcommand{\Crk}{\ensuremath{\mathbf{C}_{r,k}}}
\newcommand{\Mrk}{\ensuremath{\mathbf{M}_{r,k}}}
\newcommand{\Frk}{\ensuremath{\mathbf{f}_{r,k}}}
\newcommand{\Grk}{\ensuremath{\mathbf{g}_{r,k}}}

\newcommand{\pk}{\ensuremath{\mathbf{p}_{k}}}
\newcommand{\Vk}{\ensuremath{\mathbf{V}_{k}}}
\newcommand{\Tk}{\ensuremath{\mathbf{T}_{k}}}

\newcommand{\Krkt}{\ensuremath{\tilde{\mathbf{K}}_{r,k}}}
\newcommand{\Crkt}{\ensuremath{\tilde{\mathbf{C}}_{r,k}}}
\newcommand{\Mrkt}{\ensuremath{\tilde{\mathbf{M}}_{r,k}}}
\newcommand{\Frkt}{\ensuremath{\tilde{\mathbf{f}}_{r,k}}}
\newcommand{\Grkt}{\ensuremath{\tilde{\mathbf{g}}_{r,k}}}
\newcommand{\Vkt}{\ensuremath{\tilde{\mathbf{V}}_k}}

\begin{document}
\title{Inconsistency Removal of Reduced Bases in Parametric Model Order Reduction by Matrix Interpolation using Adaptive Sampling and Clustering} 

\author{Sebastian Resch-Schopper$^1$, Romain Rumpler$^2$, and Gerhard Müller$^1$}

\date{\today \\[2mm]
    $^1$Technical University of Munich, Chair of Structural Mechanics, Arcisstr. 21, 80333 Munich (sebastian.resch-schopper@tum.de, gerhard.mueller@tum.de)\\%
    $^2$KTH Royal Institute of Technology, MWL Laboratory for Sound and Vibration Research, Teknikringen 8, 10044 Stockholm (rumpler@kth.se)\\[2mm]
}

\begin{abstract}
Parametric model order reduction by matrix interpolation allows for efficient prediction of the behavior of dynamic systems without requiring knowledge about the underlying parametric dependency. Within this approach, reduced models are first sampled and then made consistent with each other by transforming the underlying reduced bases. Finally, the transformed reduced operators can be interpolated to predict reduced models for queried parameter points. However, the accuracy of the predicted reduced model strongly depends on the similarity of the sampled reduced bases. If the local reduced bases change significantly over the parameter space, inconsistencies are introduced in the training data for the matrix interpolation. These strong changes in the reduced bases can occur due to the model order reduction method used, a change of the system's dynamics with a change of the parameters, and mode switching and truncation. In this paper, individual approaches for removing these inconsistencies are extended and combined into one general framework to simultaneously treat multiple sources of inconsistency. For that, modal truncation is used for the reduction, an adaptive sampling of the parameter space is performed, and eventually, the parameter space is partitioned into regions in which all local reduced bases are consistent with each other. The proposed framework is applied to a cantilever Timoshenko beam and the Kelvin cell for one- to three-dimensional parameter spaces. Compared to the original version of parametric model order reduction by matrix interpolation and an existing method for inconsistency removal, the proposed framework leads to parametric reduced models with significantly smaller errors. 
\end{abstract}

\maketitle

\let\thefootnote\relax


\bigskip

\section{Introduction}%
\label{sec:intro}

Simulations of dynamical systems using the Finite Element Method (FEM) can be computationally expensive because large systems of equations have to be solved for multiple instances in time or frequency. This computational effort increases even further if the system must be evaluated for several different parameter realizations, which is required in multi-query applications like optimization or uncertainty quantification. Therefore, methods that reduce the computational effort of these applications are needed. \\
Projection-based model order reduction (MOR) can be used to accelerate the simulation for a specific realization of the parameters. These methods project the high-dimensional solution onto a lower-dimensional subspace suited to approximate the solution. To carry out this projection, the high-dimensional (full) operators are pre- and postmultiplied by a reduced basis, that spans the lower-dimensional subspace. This leads to reduced operators that describe a system that is significantly smaller but ideally still accurate. \cite{Benner2015}\\
For parametric applications, the reduced model should also maintain the parametric dependency of the high-dimensional system. For this purpose, parametric model order reduction (pMOR) can be used. These methods can be classified into global and local methods. In the former, the same reduced basis is used for the whole parameter space \cite{Benner2015}, while in the latter, some reduced quantities, such as the reduced bases \cite{Amsallem2008}, the reduced operators \cite{Amsallem2011, Panzer2010}, or the reduced transfer functions \cite{Baur2009}, are interpolated. For global methods to be efficient, the parametric dependency should exhibit a form of generalized affine representation, \textit{i.e.} a linear combination of constant matrices with scalar functions \cite{Benner2015}. This can often be achieved for dependencies involving material parameters. For geometric parameters, however, an efficient affine representation may be difficult or even impossible to obtain. \\
In this paper, we use a local pMOR method based on matrix interpolation \cite{Amsallem2011, Panzer2010}, which can treat parametric dependencies where no efficient affine representation is available. Instead, reduced models are first computed independently for a set of sample points in the parameter space. Next, the reduced operators are transformed such that they reside in the same or a similar reduced basis. Finally, the transformed reduced operators are interpolated such that a reduced system can be predicted for queried parameter points. For this final step, a similarity of the underlying bases is of utmost importance. If the transformed reduced bases are not similar, inconsistencies occur in the training data for the matrix interpolation as the reduced operators are described in different reduced bases. In turn, the resulting predicted reduced models are inaccurate. \\
In the literature, several reasons for these inconsistencies have been detected. The MOR method used \cite{Fischer2015}, the inherent change of the dynamics of a system \cite{Bazaz2015, Varona2017} and mode switching and truncation \cite{Amsallem2015} are examples of such causes for potential inconsistencies. To avoid these inconsistencies, using modal methods for generating the reduced basis \cite{Fischer2015}, performing an adaptive sampling \cite{Bazaz2015, Varona2017} and truncating the reduced bases further \cite{Amsallem2015} have been proposed. However, these remedies only cure the reason for the inconsistency they are specifically targeting. Furthermore, some of these remedies are limited in their applicability. For example, the adaptive sampling proposed in \cite{Bazaz2015, Varona2017} requires the sample points to be distributed in a regular grid. This leads to a large amount of samples added in each step of the adaptive sampling since this regular grid must be maintained. Furthermore, treating inconsistencies due to mode switching and truncation by removing parts of the reduced bases deteriorates the accuracy of the sampled reduced systems. \\
In this paper, we thus combine multiple methods for treating inconsistencies in the reduced bases into a general framework: For MOR, a modal method is used. Furthermore, the adaptive sampling is extended to randomly sampled points by identifying neighboring samples via a triangulation. This allows placing only one new sample point at a time, as there is no structure of the sample points that must be maintained. Finally, the inconsistency due to mode switching and truncation is treated by partitioning the parameter space into several regions such that all local reduced bases within one region are consistent with each other. This avoids truncating the reduced bases further and, thus, does not deteriorate the accuracy of the sampled reduced systems. The idea of partitioning the parameter space and constructing several local parametric reduced order models (pROMs) has already been used in the literature for other pMOR approaches \cite{Amsallem2012, Eftang2012, Haasdonk2011, Peherstorfer2014} and is here implemented with the other methods combined in one general framework. \\
The remaining paper is structured as follows: Section \ref{sec:Background} introduces the systems that are investigated and the concept of (parametric) model order reduction. In section \ref{sec:Inconsistencies}, reasons and remedies for inconsistencies in the reduced bases are explained. Section \ref{sec:AdaptiveSampling} presents the proposed framework for removing inconsistencies. Finally, the results of this approach are shown in section \ref{sec:Results}, followed by a conclusion in section \ref{sec:conclusions}.

\section{Theoretical Background}%
\label{sec:Background}

\subsection{Problem Definition}

We consider parameter-dependent, linear time-inva\-riant dynamical systems with single input and single output (SISO) in second-order form in frequency domain: 
\begin{equation}
   \Sigma(\p): \begin{cases}
		s^2\M(\p) \x(\p,s) + s\C(\p) \x(\p,s) + \K(\p) \x(\p,s) &= \F(\p) u(s), \\
		\hfill y(\p,s) &= \G(\p)\mathbf{x}(\p,s),
		\end{cases}
\end{equation}
with parameter-dependent mass, damping and stiffness matrices $\mathbf{M}, \, \mathbf{C}, \, \mathbf{K}: \mathbb{R}^d \rightarrow \mathbb{R}^{n \times n}$, degrees of freedom $\x \in \mathbb{R}^n$, complex frequency $s \in \mathbb{C}$, input $u \in \R$ and $\F: \R^d \rightarrow \R^{n \times 1}$ and output $y \in \R$ and $\G: \R^d \rightarrow \R^{1 \times n}$. $\M, \C, \K, \F$ and $\G$ are assumed to depend on $d$ parameters $\p = [p_1, p_2, \dots , p_d] \in \R^d$; $n$ denotes the number of degrees of freedom. For the sake of simplicity, we restrict ourselves to single input single output (SISO) systems, but the concepts presented in the following may be applied to multiple input multiple output (MIMO) systems as well.\\

\subsection{Projection-based Model Order Reduction}
\label{sec:MOR}

In the following we assume that the full order model (FOM) is sampled at $K$ parameter points $\p_k$ with $k=1, \dots, K$. To reduce the computational effort for solving a FOM at one parameter point $\pk$, the full system's response $\x(s)$ can be projected onto a lower dimensional subspace $\mathcal{V}_k$, which is spanned by the columns of a reduced basis $\V_k \in \R^{n \times r}$. When using a Galerkin-projection, this results in the following reduced model of size $r \ll n$:
\begin{equation}
    \Sigma_r(\pk): \begin{cases}
		s^2\Mrk \x_r(s) + s\Crk \x_r(s) + \Krk \xr(s) &= \Frk u(s), \\
		\hfill y_r(s) &= \Grk \xr(s),
		\end{cases} \label{eq:ReducedSystem}
\end{equation}
where 
\begin{equation}
\begin{aligned}
    \Mrk &= \Vk^\intercal \M(\pk) \Vk, \qquad &\Crk &= \Vk^\intercal \C(\pk) \Vk, \qquad \Krk = \Vk^\intercal \K(\pk) \Vk \\
    \Frk &= \Vk^\intercal \F(\pk), \qquad &\Grk &= \G(\pk) \Vk.
    \label{eq:Projection}
\end{aligned}
\end{equation}
By reformulating system \cref{eq:ReducedSystem}, the transfer function of the reduced model at the sample point $\p_k$ is obtained as
\begin{equation}
	H_{r,k}(s) = \Grk \left( s^2 \Mrk + s \Crk + \Krk \right)^{-1} \Frk.
\end{equation}
For generating the reduced basis $\Vk$, various methods can be applied. Popular choices are modal methods \cite{Tiso2021}, moment matching \cite{Benner2021a} or Proper Orthogonal Decomposition (POD) \cite{Antoulas2005}. In this work, a modal method is used like in several contributions on pMOR by matrix interpolation \cite{Panzer2010, Fischer2014, Fischer2015, Mencik2020}.\\
 In modal truncation (MT), the reduced basis $\mathbf{V}$ consists of selected eigenmodes, which are the eigenvectors $\mathbf{\Phi}$  of the eigenvalue problem associated with the undamped system: \cite{Tiso2021}
\begin{equation}
	(\omega^2 \mathbf{M} + \mathbf{K}) \mathbf{\Phi} = \mathbf{0},
\end{equation}
where $\omega$ denotes the angular eigenfrequencies of the system. One possibility to build the reduced basis is to select the eigenmodes associated with the $r$ smallest eigenfrequencies, which results in a reduced model targeting accuracy for low frequencies. Alternatively, the $r$ eigenmodes that contribute most to the input-output-behavior based on  a dominance criterion may be selected \cite{Rommes2008}.

 \subsection{Parametric Model Order Reduction by Matrix Interpolation}

 To maintain the parametric dependency in the reduced model, parametric model order reduction (pMOR) can be applied. In this work, we assume that there is no access to an (efficient) affine representation of the parametric dependency. Therefore, pMOR by matrix interpolation \cite{Panzer2010} is used. This approach starts by computing reduced models for $K$ sample points in the parameter space, as outlined in \cref{sec:MOR}. The objective is ultimately to interpolate the reduced operators. However, since each of the sampled models is reduced independently, they are described in different reduced bases. Thus, an interpolation of the reduced operators does not make sense yet. Therefore, a reference basis is required first, to which all sampled reduced bases shall be transformed. For this purpose, the Singular Value Decomposition (SVD) is performed on the library of collected local reduced bases $\Vk$: 
 \begin{equation}
     \mathbf{U} \pmb{\Sigma} \mathbf{Z} = [\V_1, \V_2, \dots, \V_K].
     \label{eq:Vall}
 \end{equation}
 The reference basis $\mathbf{R}$ is then obtained as the first $r$ columns of $\mathbf{U}$. Next, a linear transformation is applied to each local reduced basis:
  \begin{equation}
     \Vkt = \Vk \Tk,
     \label{eq:TransformationBasis}
 \end{equation}
 where $\Vkt$ is the transformed reduced basis and $\Tk$ the transformation matrix. The transformed reduced basis should be as similar as possible to $\mathbf{R}$. This is the case if $\Vkt^\intercal \mathbf{R} = \mathbf{I}$, where $\mathbf{I}$ is the identity. Therefore, $\Tk$ is defined as 
 \begin{equation}
    \Tk = (\mathbf{R}^\intercal \Vk)^{-1}. \label{eq:Tk}
\end{equation}
The transformed reduced operators are then obtained by applying this transformation:
\begin{equation}
\begin{aligned}
    \Mrkt &= \Tk^\intercal \Mrk \Tk, \qquad &\Crkt &= \Tk^\intercal \Crk \Tk, \qquad \Krkt = \Tk^\intercal \Krk \Tk    \\
    \Frkt &= \Tk^\intercal \Frk, \qquad &\Grkt &= \Grk \Tk.
    \label{eq:TransfOp}
\end{aligned}
\end{equation}
Finally, the reduced operators can be interpolated. This can either be done entry-wise or at once for the complete matrix. For the interpolation itself, any interpolation or regression method can be used. In order to guarantee that the predicted reduced operators fulfill important properties such as symmetry and positive definiteness, the interpolation can be carried out on the tangent space of the manifold the reduced operators belong to \cite{Amsallem2011}. As an alternative, the Cholesky decomposition can be computed, which is a factorization for symmetric positive definite matrices. Instead of interpolating the entries of the mass, damping and stiffness matrix, the entries of their Cholesky factors are then interpolated \cite{Amsallem2015}. \\

\section{Inconsistencies in Reduced Bases}%
\label{sec:Inconsistencies}

In the final step of pMOR by matrix interpolation, the entries of the transformed reduced operators are interpolated. According to \cref{eq:TransfOp}, the dependency to be interpolated does not only consist of the parametric dependency of the full operators but also on the parametric dependency of the reduced bases:
\begin{equation}
    \tilde{\mathbf{K}}_r(\p) = \tilde{\V}^{\trans}(\p) \K(\p) \tilde{\V}(\p).
\end{equation}
In the general case, $\tilde{\V}(\p)$ is not constant but may differ for different points in the parameter space. The transformation in \cref{eq:TransformationBasis} does not necessarily remove this dependency since it is a linear transformation that does not change the subspace that $\Vk$ spans. Consequently, strong changes in the subspaces that the sampled reduced bases $\Vk$ span, can affect the total dependency $\tilde{\mathbf{K}}_r(\p)$ significantly, which introduces inconsistencies in the training data for the matrix interpolation. In order to improve the accuracy of the interpolation, these inconsistencies should be minimized. Therefore, detecting them and understanding their origin is important.

\subsection{Principal Subspace Angles}

Inconsistencies are large changes in the sampled reduced bases. Thus, a measure of the distance between the subspaces spanned by the reduced bases is required to detect them. The difference between two orthonormal bases $\V_i$ and $\V_j$ can be assessed by the principal angles between these subspaces. For computing them, the SVD of the product $\V_i^{\trans} \V_j$ is required first: \cite{Amsallem2015}
\begin{equation}
    \mathbf{W}_i \pmb{\Sigma} \mathbf{Z}_j = \V_i^{\trans} \V_j,
    \label{eq:SubspaceAngles}
\end{equation}
where $\mathbf{W}_i = [\mathbf{w}_{i,1}, \dots, \mathbf{w}_{i,r}]$, $\mathbf{Z}_j = [\mathbf{z}_{j,1}, \dots, \mathbf{z}_{j,r}]$, and $\pmb{\Sigma} = \mathrm{diag}(\sigma_1, \dots, \sigma_r)$. From this, the principal angles are obtained as \cite{Amsallem2015}
\begin{equation}
    \theta_l = \arccos(\sigma_l), \quad l = 1, \dots r.
\end{equation}
The principal angle $\theta_l$ refers to the angle between the vectors $\V_i \mathbf{w}_{i,l}$ and $\V_j \mathbf{z}_{j,l}$. We assume that the singular values are ordered decreasingly and are bounded by the interval $[0, 1]$, thus it holds that $0~\leq~\theta_1~\leq~\dots~\leq~\theta_r~\leq~90^\circ$.

\subsection{Sources of Inconsistencies and Possible Remedies}

Since the reduced basis is generated independently for each sampled parameter point, there are several reasons that can lead to great differences in the reduced basis.

\subsubsection{Model Order Reduction Method}

\begin{figure}[htbp]
	\centering
	\includegraphics[width=\textwidth]{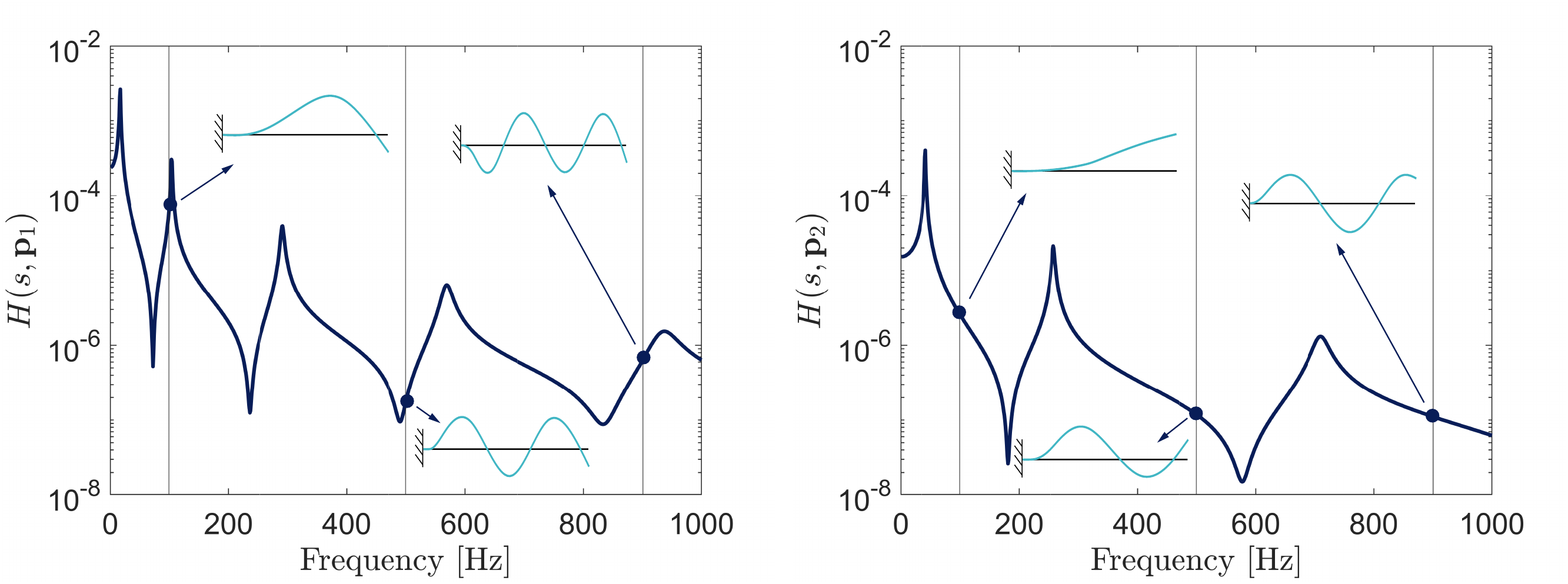}
	\caption{Snapshots generated by Proper Orthogonal Decomposition for a 2D cantilever beam for two different realizations of the parameters. For each snapshot frequency, the physical deformation of the beam is shwon in the small sketch next to it.}
	\label{fig:POD_Beam_Final}
\end{figure}

The MOR method used has a major impact on the similarity between the sampled reduced bases. In order to minimize changes in the sampled reduced bases $\Vk$, the MOR method used should generate the same or similar basis vectors for different realizations of the parameters. Snapshot-based MOR methods like POD, for example, do not fulfill this property since the collected basis vectors depend on the frequencies used for generating the snapshots. We illustrate this for a 2D cantilever beam with a transversal harmonic excitation in \cref{fig:POD_Beam_Final}. If a parameter of the beam is changed, the dynamic response at a given frequency may be greatly affected. Contrary to POD, modal methods identify the eigenmodes of the structure without requiring the user to input frequencies for generating the reduced basis. In the example of the cantilever beam, modal truncation would generate the same basis vectors for the different realizations of the parameters, being consistently based on the bending modes with the smallest eigenfrequencies. Similar findings are reported in \cite{Fischer2015} for the simulation of moving loads. Therein, it was reported that the dependency of the transformed reduced operators was smoother when using Component Mode Synthesis compared to Rational Interpolation. 

\subsubsection{Change of the System Dynamics} \label{sec:ChangeSystemDynamics}

Furthermore, inconsistencies may be caused due to a change of the underlying system dynamics with a change of the parameters. If regions where the system dynamics change strongly are sampled poorly, the sampled reduced bases will also differ strongly. This source of inconsistency is independent of the MOR method used since the change is inherent to the system. This issue has been investigated and addressed in \cite{Bazaz2015, Varona2017}. The authors of these contributions propose to use adaptive sampling to solve the problem: By computing the principal subspace angles between neighboring samples, the difference between two reduced bases can be quantified. If the largest subspace angle is above a user-defined threshold, new samples are placed in this region. This procedure is repeated until the largest subspace angle is below the threshold for all neighboring samples. Consequently, regions in the parameter space where the reduced bases change strongly are sampled with a fine resolution so that the change of the reduced bases can be captured well.

\subsubsection{Mode Switching and Truncation}

\begin{figure}[htbp]
	\centering
	\begin{subfigure}[b]{0.49\textwidth}
		\centering
		\includegraphics[width=0.8\textwidth]{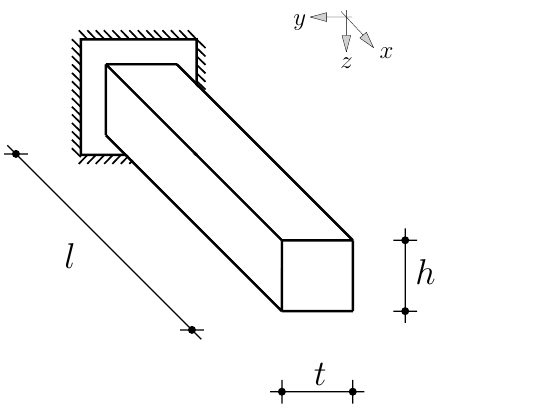}
	\end{subfigure}
	\hfill
	\begin{subfigure}[b]{0.49\textwidth}
		\centering
		\includegraphics[width=\textwidth]{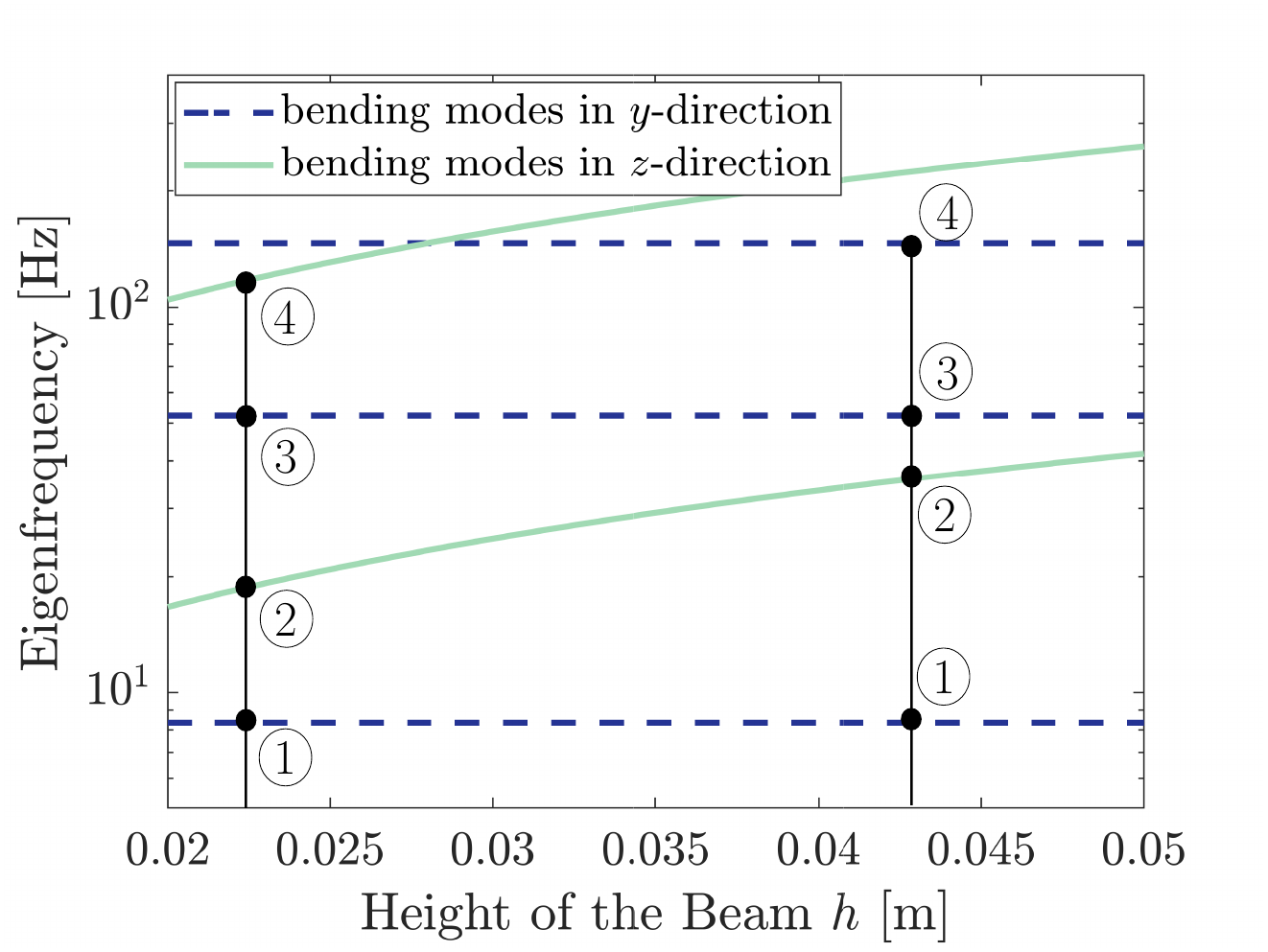}
	\end{subfigure}
	\caption{Eigenfrequencies of a cantilever bending beam for changing height of the beam.}
	\label{fig:TruncationSwitch}
\end{figure}

Finally, an inconsistency may occur due to mode switching and truncation: With a change of the parameters, some eigenmodes may switch their position. In modal MOR methods, only a limited amount of eigenmodes is kept for generating the basis, thus, parts of the switched modes may be removed from the reduced basis for some realizations of the parameters \cite{Amsallem2015}. This is illustrated in \cref{fig:TruncationSwitch} for a 3D cantilever beam where the beam's height is used as a parameter. 
While an increase in the height of the beam leads to higher eigenfrequencies of the bending modes in the $z$-direction, the eigenfrequencies of the bending modes in the $y$-direction are unaffected. Consequently, the second bending mode in the $z$-direction and the third bending mode in the $y$-direction swap positions at a specific value of the height. If only the first four eigenmodes are kept for generating the reduced model, the reduced basis would consist of two bending modes in the $z$-direction and two bending modes in the $y$-direction for smaller heights, whereas it would contain one bending mode in the $z$-direction and three bending modes in the $y$-direction for larger heights. In the case of modal truncation, the subspaces spanned by these two inconsistent bases are even orthogonal. This, in turn, introduces numerical errors in the transformation because the product $\mathbf{R}^\intercal \Vk$ and thus the transformation matrix $\Tk$ become singular. \\
This problem has been identified in \cite{Amsallem2015}. In order to remove this inconsistency, the authors propose to compute the subspace angles between all sampled reduced bases and the reference basis. For each reduced basis, the number of angles smaller than $45^\circ$ is denoted with $l_k$ and stored. Furthermore, the matrix $\mathbf{W}_k$ from the computation of the subspace angles in \cref{eq:SubspaceAngles} is stored for each reduced basis $\Vk$. Finally, the smallest $l_k$ is identified, and all reduced bases are truncated to the first $l_k$ columns of $\Vk \mathbf{W}_k$. \\
Due to this further truncation, the accuracy of the sampled reduced models deteriorates. Furthermore, in case the problem of mode switching and truncation occurs multiple times in the parameter space, the method will not remove all inconsistencies. This is because the parts of one reduced basis that are consistent with the reference basis are not necessarily the same for another reduced basis. 
\section{Adaptive Sampling and Clustering for Treatment of Inconsistencies}
\label{sec:AdaptiveSampling}

In the following, the independent remedies for the individual sources of inconsistencies are partially adapted and combined into a general workflow. For this, an adaptive sampling is performed first to minimize the inconsistencies stemming from a change of the underlying system dynamics and to identify inconsistencies from mode switching and truncation. In the sampling, the reduced models are obtained by modal truncation. Afterwards, the samples are clustered such that the reduced bases of all samples within a cluster are consistent with each other. Based on this clustering, the parameter space is partitioned into consistent regions. Due to this, the local reduced models do not lose accuracy and multiple mode switching can be treated. Finally, the individual consistent regions can be filled with further samples if needed. In the following, the individual steps are explained.

\subsection{Adaptive Sampling}

\begin{figure}[thb]
	\centering
	\includegraphics[width=0.8\textwidth]{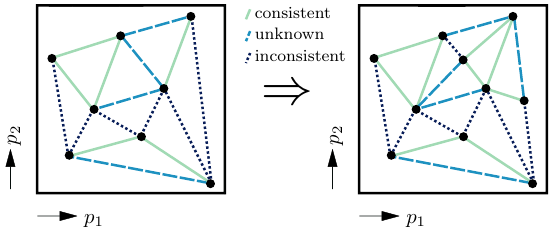}
	\caption{Schematic representation of adaptive sampling: Consistency and samples before and after one iteration of adaptive sampling.}
	\label{fig:AdaptiveSampling}
\end{figure}

The adaptive sampling is based on the subspace angles between neighboring sample points as proposed in \cite{Bazaz2015, Varona2017} and explained in \cref{sec:ChangeSystemDynamics}. In order to extend this approach to randomly distributed sample points, the Delaunay triangulation \cite{Delaunay1934} is used to identify neighboring samples. The criterion, where to place new samples, is based on an upper and a lower threshold on the subspace angles so that two neighboring samples are considered consistent if their subspace angle is smaller than the lower threshold and inconsistent if it is larger than the upper threshold. If the largest subspace angle is between the two thresholds, the consistency of the two samples is not known yet and needs to be determined in the course of the adaptive sampling. Furthermore, thresholds for the Euclidian distance between neighboring samples in the normalized parameter space are used to avoid large unsampled regions and to prevent the algorithm from getting stuck by placing many samples in small parts of the parameter space. \\
The adaptive sampling procedure is outlined in \cref{alg:AdaptiveSampling} and shown schematically in \cref{fig:AdaptiveSampling}. As inputs, initial sample points $\mathbf{P} = [\p_1, \dots , \p_K]$ are required, which can be distributed regularly or randomly. In the tested examples, a regular distribution of the initial samples with $2^d$ samples is used so that the initial samples are only the corner points of the parameter space. This allows the algorithm to find the "optimal" distribution of samples in the parameter space based on the criteria for the adaptive sampling. Furthermore, the reduced systems computed at the sample points are required. Finally, scalar values for the upper and the lower threshold for the subspace angles $\theta_{uT}$ and $\theta_{lT}$ and the thresholds $d_{uT}$, $d_{lT}$ and $d_N$ for the distance between samples must be defined. 
\begin{algorithm}[thb]
	\caption{Adaptive Sampling}\label{alg:AdaptiveSampling}
	\begin{algorithmic}[1]
		\Require Initial sample points $\mathbf{P}$, sampled reduced systems $\Sigma_{r}$, upper and lower threshold for the subspace angles $\theta_{uT}$ and $\theta_{lT}$, thresholds for the distance between neighboring samples $d_{uT}$, $d_{lT}$ and $d_N$
		\Ensure sample points $\mathbf{P}$, reduced systems $\Sigma_{r}$
		\State get edges of Delaunay triangulation of $\mathbf{P}$ \label{lin:Delaunay}
		\For{$i = 1$ to number of edges}
		\State $d_i \gets$ length of edge $i$
		\State $\theta_i \gets$ largest principal subspace angle between samples of edge $i$
		\State $c_i \gets$ getConsistency($d_i, \, \theta_i, \, \theta_{uT}, \, \theta_{lT}, \, d_{lT}$)
		\EndFor
		
		\While{any $c_i = 0$ or any $d_i$ of inconsistent neighbors $> d_{uT}$}
		
		\If{any $c_i = 0$} \label{lin:c_i}
		\State\label{lin:test}{get all pairs of neighbors for which $c_i = 0$}
		\State find the two neighbors that have the largest $\theta_i$ and for which the distance \Statex \hspace{1.5cm} between the mid-point of this edge and all other samples is larger than \Statex \hspace{1.5cm} $d_{N}$
		\State place new sample point $\mathbf{p}_\theta$ between them and compute reduced system \Statex \hspace{1.5cm} $\Sigma_r(\mathbf{p}_\theta)$
		\State $\mathbf{P} \gets [\mathbf{P}, \mathbf{p}_\theta]$
		\State $\Sigma_r \gets [\Sigma_r, \Sigma_r(\mathbf{p}_\theta)]$
		\EndIf
		
		\If{any $d_i$ of inconsistent neighbors $> d_{uT}$} \label{lin:d_i}
        \State get all pairs of neighbors for which $c_i = 0$
		\State find the two neighbors that have the largest $d_i$ and for which the distance \Statex \hspace{1.5cm} between the mid-point of this edge and all other samples is larger than \Statex \hspace{1.5cm} $d_{N}$
		\State place new sample point $\mathbf{p}_d$ between them and compute reduced system \Statex \hspace{1.5cm} $\Sigma_r(\mathbf{p}_d)$
		\State $\mathbf{P} \gets [\mathbf{P}, \mathbf{p}_d]$
		\State $\Sigma_r \gets [\Sigma_r, \Sigma_r(\mathbf{p}_d) ]$
		\EndIf
		
		\State get edges of Delaunay triangulation
		\For{$i = 1$ to number of edges}
		\State $d_i \gets$ length of edge $i$
		\State $\theta_i \gets$ largest principal subspace angle between samples of edge $i$
		\State $c_i \gets$ getConsistency($d_i, \, \theta_i, \, \theta_{uT}, \, \theta_{lT}, \, d_{lT}$)
		\EndFor
		
		\If{a $\mathbf{p}_\theta$ was computed in this iteration}
		\State $c_i \gets$ update consistency based on change of subspace angle
		\EndIf
		\EndWhile
	\end{algorithmic}
\end{algorithm}
\begin{algorithm}[thb]
	\caption{getConsistency}\label{alg:getConsistency}
	\begin{algorithmic}[1]
		\Require Distance $d$ between the two neighboring samples, largest subspace angle $\theta$ between the two neighboring samples, upper and lower threshold for the subspace angles $\theta_{uT}$ and $\theta_{lT}$, lower threshold for the distance between neighboring samples $d_{lT}$
		\Ensure Consistency $c$ of the two neighboring samples
		\If{$\theta \leq \theta_{lT}$}
		\State $c \gets 1$    \Comment{Samples are consistent}
		\ElsIf{$\theta \geq \theta_{uT}$ or $d < d_{lT}$}
		\State $c \gets -1$   \Comment{Samples are inconsistent}
		\Else
		\State $c \gets 0$    \Comment{Consistency of samples is unknown yet}
		\EndIf
	\end{algorithmic}
\end{algorithm}

In the first step of the algorithm, the neighbors of all samples are found by computing the Delaunay triangulation in line \ref{lin:Delaunay}. This returns a triangular mesh of the sample points so that a list of pairwise neighboring samples is obtained. In MATLAB, a Delaunay triangulation may be performed in 2D or 3D with the function \verb+delaunayTriangulation+ and in higher dimensions with \verb+delaunayn+. Note that the parameter space should be normalized to the interval $[0, \, 1]$ to avoid the Delaunay triangulation of the initial samples being distorted due to different magnitudes of the parameters. Afterwards, the distance $d_i$ and the largest subspace angle $\theta_i$ between all neighboring samples in the normalized parameter space are computed. With this information, a flag for the consistency $c_i$ can be evaluated according to the rules stated before and shown in \cref{alg:getConsistency}. Note that the only physical criterion for inconsistency is the subspace angle but the criterion for the distance is used to avoid placing many samples in small regions of the parameter space. If none of the above criteria are fulfilled, the consistency is unknown yet and $c_i$ is set to 0.  \\
Afterwards, the adaptive sampling is performed based on the two criteria stated in lines \ref{lin:c_i} and \ref{lin:d_i}: If there are two neighboring samples with unknown consistency, a new sample is placed between them. For that criterion, the edge with the largest subspace angle is selected so that the neighbors with the highest difference in the reduced bases are checked first. Furthermore, a sample is placed between neighboring inconsistent samples whose distance is larger than the threshold $d_{uT}$. For both criteria, also the distance between the mid-point of this edge and all other samples is checked and compared to $d_N$, a threshold for the distance to the neighbors. If the mid-point is closer to a neighboring sample than $d_N$, no new sample is placed at the mid-point of this edge. This prevents the algorithm from placing many sample points in small regions of the parameter space and is only relevant for three- or higher-dimensional parameter spaces, where skew edges can occur. Due to these, there might be edges whose lengths are larger than $d_{lT}$ but whose mid-points are close to already existing sample points. Since new samples have been added, the triangulation and consequently $d_i$, $\theta_i$, and $c_i$ must be updated. This procedure is also illustrated in \cref{fig:AdaptiveSampling}: A new sample is placed in the middle of the edge with unknown consistency with the largest subspace angle. This could, for example, be the edge in the middle upper area. Afterwards, the neighbors of the new sample are identified and it turns out that the new sample is consistent with its neighbors to the right, inconsistent with the neighbor at the top left and the consistency to the neighbor to the bottom left is unknown yet. Furthermore, one sample is placed in the middle of the right edge, which is the longest inconsistent edge. Again, the new neighbors are identified and the consistencies are determined. 

Finally, it is important to point out a further criterion for inconsistency that can occur when a sample $\mathbf{p}_\theta$ is added due to unknown consistency: So far, two neighboring samples are only considered inconsistent if their largest subspace angle is above the threshold $\theta_{uT}$ or their distance is below $d_{lT}$. However, the largest subspace angle may also hardly decrease when placing a new sample between two neighboring samples with unknown consistency. This could cause the algorithm to get stuck. Therefore, the change of the subspace angles is checked when a sample $\mathbf{p}_\theta$ is added. If the subspace angle did not decrease by more than a certain amount, for example, $10\%$, the two samples are tagged as inconsistent.

\subsection{Clustering}

\begin{figure}[htb]
	\centering
	\includegraphics[width=0.8\textwidth]{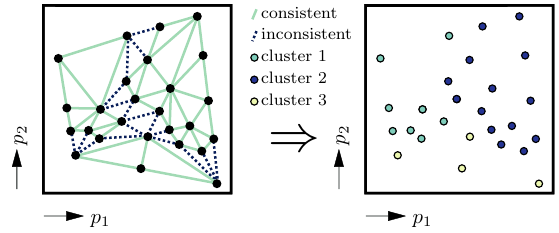}
	\caption{Schematic representation of clustering: Consistency of samples after adaptive sampling and clustering of samples.}
	\label{fig:Clustering}
\end{figure}

When the consistencies are known for all neighboring samples, the samples can be clustered. The steps of this procedure are outlined in \cref{alg:Clustering} and the result is shown schematically in \cref{fig:Clustering}. In this procedure, a loop over all edges of the Delaunay triangulation is performed and for each edge, the two samples of this edge are assigned to the same class if they are consistent. Since the edges are checked in an arbitrary order, two samples of the same class can be assigned to different classes during the algorithm. Therefore, this is detected in line \ref{lin:TwoClasses} and the superfluous class is deleted in line \ref{lin:DeleteClass}. Finally, all samples that are not consistent with any other sample are assigned to a separate class.
\begin{algorithm}[thb]
\caption{Clustering}\label{alg:Clustering}
\begin{algorithmic}[1]
\Require Edges of the Delaunay triangulation, consistency $c$ of all neighboring samples 
\Ensure Class of each sample
\For{$i=1$ to number of edges}
    \If{$c_i = 1$}
        \If{none of the two samples is assigned to a class yet}
            \State assign both samples to a new class
        \ElsIf{only one of the two samples is already assigned to a class}
            \State assign the unassigned sample to the same class as the other sample
        \ElsIf{both samples are already assigned, but to different classes} \label{lin:TwoClasses}
            \State delete one of the two class labels and assign all samples of the deleted \Statex \hspace{2cm} class to the not deleted class \label{lin:DeleteClass}
        \EndIf
    \EndIf
\EndFor

\For{$i=1$ to number of samples}
    \If{sample $i$ is not assigned to a class}
        \State assign sample to a new class
    \EndIf
\EndFor
\end{algorithmic}
\end{algorithm}

\subsection{Fill Clusters}

\begin{figure}[htb]
	\centering
	\includegraphics[width=0.8\textwidth]{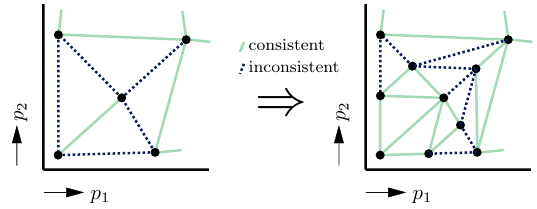}
	\caption{Schematic representation of finding the borders of a region: Consistency of samples before and after more samples are added to identify the border between two inconsistent regions.}
	\label{fig:FindBorders}
\end{figure}

The adaptive sampling is solely based on the subspace angles and the distance between neighboring samples. Therefore, some clusters may consist of only one or a few samples. In this case, performing an interpolation of the reduced operators might not make sense or can, depending on the method used, even be impossible - for example, when using polynomial interpolation or regression with a required order. Therefore, these regions must be filled with further samples. This can be done by successively placing further samples between the samples of this region with the largest distance until the required number of samples is achieved. These steps are also outlined in \cref{alg:FillRegion}. However, a special case needs further treatment: when the cluster consists of fewer or as many samples as the dimension of the parameter space or if the samples of the cluster do not span the complete space, further samples would not be distributed in all dimensions but only in the dimensions spanned by the already existing sample points. For the case of a single sample within a cluster, the approach outlined in \cref{alg:FillRegion} would not work at all. Therefore, the borders of this region need to be found first in all dimensions, which is shown in \cref{alg:FindBorders} and \cref{fig:FindBorders}. The inputs of this function are the samples of this cluster, the edges of the Delaunay triangulation and the consistency flag of all neighboring samples. If \cref{alg:FindBorders} is executed until the borders are found, this might lead to an unnecessarily dense sampling in a small region of the parameter space. 
\begin{algorithm}[tbh]
	\caption{Fill Region}\label{alg:FillRegion}
	\begin{algorithmic}[1]
		\Require Samples $\mathbf{P}$ of one cluster, edges of the Delaunay triangulation, consistency $c$ of all neighboring samples, minimum required amount of samples $m$
		\Ensure Further samples of this cluster
		\While{number of samples of this cluster $< m$}
		\State place a new sample point $\mathbf{p}$ between the two samples of this cluster with the  
		\Statex \hspace{1cm} largest distance
		\State $\mathbf{P} \gets [\mathbf{P}, \mathbf{p}]$
		\State $\Sigma_r \gets [\Sigma_r, \Sigma_r(\mathbf{p}) ]$
		\State get edges of Delaunay triangulation of samples of this cluster
		\EndWhile
	\end{algorithmic}
\end{algorithm}
Therefore, a fixed number of times the algorithm is called can be set, and if not enough samples of this cluster are found, the complete cluster is deleted and no more samples are placed in this region. This prevents the algorithm from getting stuck by placing many sample points in a small region at the cost of less accuracy in this small region, as the local pROM of a neighboring region will be used to make predictions in this region.  \\
Finally, there might still be large unsampled regions if the consistent samples of a cluster are all far apart from each other. In order to prevent this, a maximum distance between consistent samples $d_C$ can be set and samples are added in the middle of all edges between consistent samples that are further apart from each other than $d_C$.

\begin{algorithm}[htb]
\caption{Find Borders}\label{alg:FindBorders}
\begin{algorithmic}[1]
\Require Samples $\mathbf{P}$ of one cluster, edges of the Delaunay triangulation, consistency $c$ of all neighboring samples
\Ensure Further samples of this cluster
\For{$i=1$ to number of samples in this region}
    \State get the neighboring samples of sample $i$
    \For{$j=1$ to number of neighboring samples that are not part of this cluster}
        \State place new sample point $\mathbf{p}$ between samples $i$ and $j$ and compute reduced  
        \Statex \hspace{1.5cm} system $\Sigma_r(\mathbf{p})$
        \State $\mathbf{P} \gets [\mathbf{P}, \mathbf{p}]$
        \State $\Sigma_r \gets [\Sigma_r, \Sigma_r(\mathbf{p}) ]$
        \State get edges of Delaunay triangulation
        \State assign the new sample to the corresponding class based on the largest  \Statex \hspace{1.5cm} principal subspace angles with the neighboring samples
    \EndFor
\EndFor
\end{algorithmic}
\end{algorithm}

\subsection{Classification}

Now that the reduced bases of all samples within a cluster are consistent with each other and each cluster contains sufficiently many samples, the remaining steps of pMOR by matrix interpolation as outlined in \cref{eq:Vall} - \cref{eq:TransfOp} can be performed for each cluster separately. Interpolating the reduced operators of each cluster leads to several local parametric reduced order models (pROMs), that are only valid within a specific region. To find out which local model to use for a queried parameter point, a classification model is required that can predict the corresponding pROM for a queried parameter point. For this purpose, any classification method can be used. In the following, we use \textit{k-nearest-neighbor} classification. As the name suggests, this method predicts the class of a queried parameter point based on the class of the $k$ nearest neighbors. Hence, its computational effort is low and it is well suited in combination with the proposed approach because the distance between inconsistent regions is one of the hyperparameters of the method. \\
Nevertheless, there is no guarantee that for each queried point, the correct corresponding pROM is found. As long as the reduced basis of one region allows for a good approximation also in the neighboring region, this is not a problem. However, if the reduced basis of one region lacks important information required in the neighboring region, the pROM will be inaccurate. A remedy to this problem could be to compute the number of inconsistent basis vectors between two neighboring regions. If this exceeds a user-defined value, the pROM at the border of the regions should not be trusted. Instead, a global basis could be generated from the local bases of the neighboring regions by concatenating these local bases. This approach would be less efficient as it requires the full model for the prediction at a queried parameter point to perform the projection into this global basis, but it would be more accurate.

\section{Results}
\label{sec:Results}

The proposed framework of adaptive sampling and clustering is tested on two different structures with varying parameter space dimensions: a cantilever beam with one-dimensional parameter space and the Kelvin cell structural topology once with a two- and once with a three-dimensional parameter space. All computations are performed on an Intel Xeon E5-2660v3 processor with 2.6 GHz using MATLAB R2023a. The classification of the parameter space is performed using \textit{k-nearest-neighbor} classification with the default settings of MATLAB's function \verb+fitcknn+. For the one-dimensional parameter space, the reduced operators are interpolated using spline interpolation with MATLAB's function \verb+spline+. For the higher-dimensional parameter spaces, 
ridge regression with all terms up to order 3 in each variable and a ridge parameter $\lambda=10^{-5}$ was used, which proved to be good choices.

\begin{figure*}[htb]
	\centering
	\begin{subfigure}[b]{0.49\textwidth}
		\centering
		\includegraphics[width=\textwidth]{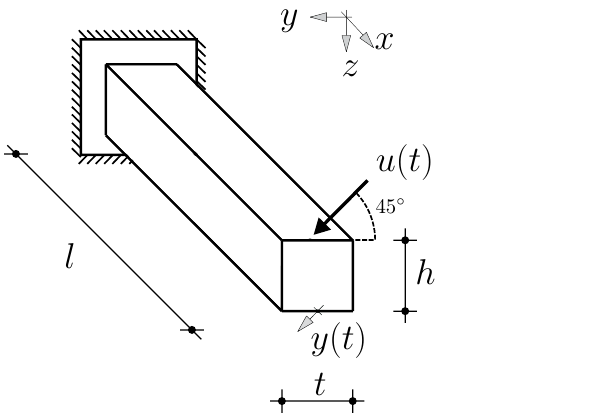}
		\caption{3D cantilever beam.}  
	   \label{fig:TimoshenkoBeam}
	\end{subfigure}
	\hfill
	\begin{subfigure}[b]{0.49\textwidth}  
		\centering 
		\includegraphics[width=\textwidth]{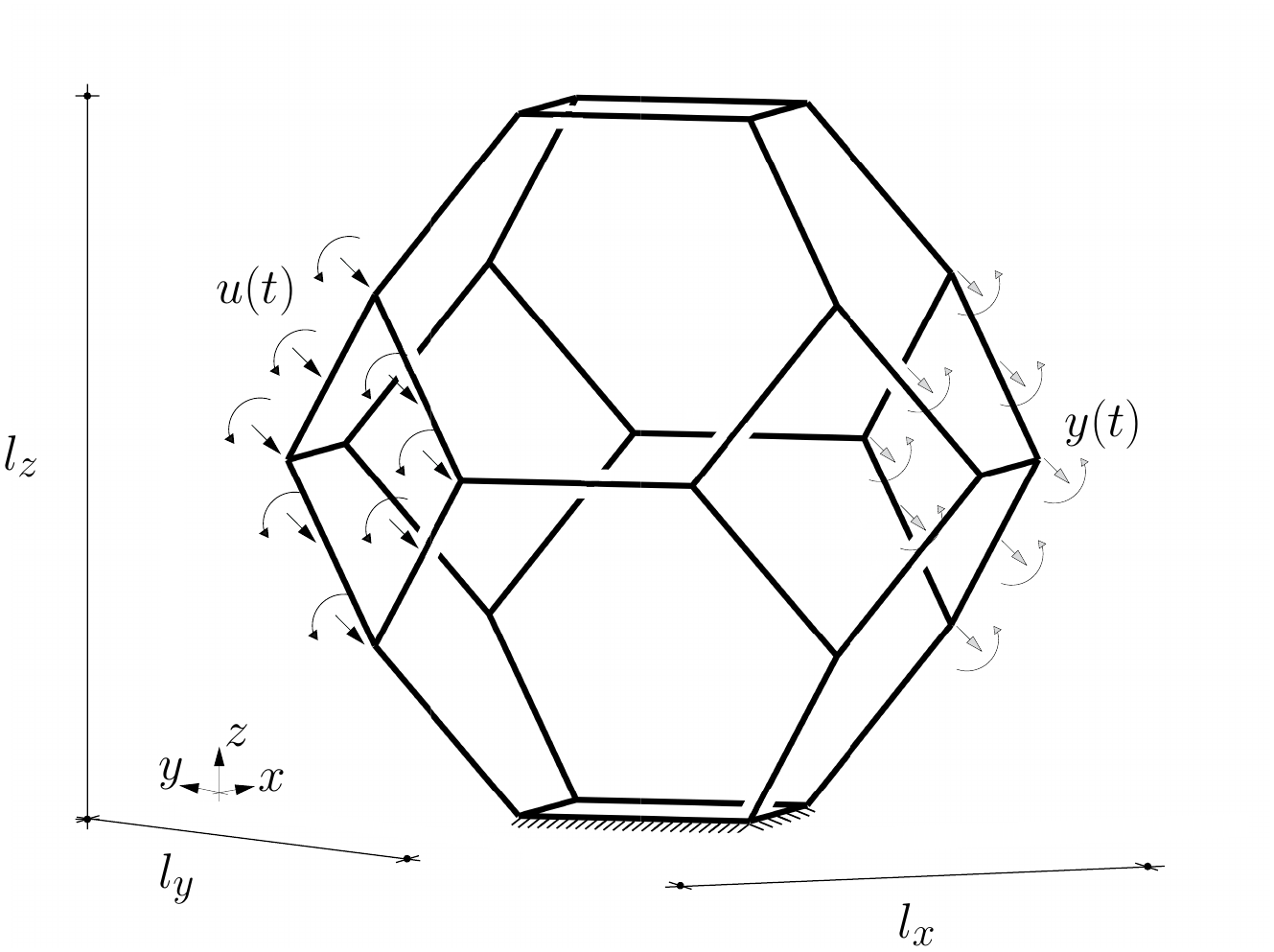}
	   \caption{Kelvin Cell.}
	   \label{fig:KelvinCell}
	\end{subfigure}
	\caption{Structures used for testing the proposed framework.}  
\end{figure*}

\subsection{Timoshenko Beam \textendash{} 1 Parameter}

First, the proposed approach is tested on an academic example, namely a 3D cantilever Timoshenko beam \cite{Panzer2009} depicted in \cref{fig:TimoshenkoBeam}. The height of the beam is chosen as parameter and varied between $2 \, \mathrm{cm}$ and $5 \, \mathrm{cm}$. The remaining parameters of the beam are chosen as shown in \cref{tab:ParametersBeam}. The beam is clamped on one end and harmonically excited at the tip with a force that is tilted in the $y$-$z$-plane in order to excite bending modes both in $y$- and in $z$-direction. Since only the height but not the width of the beam is changed, mode switching of the bending modes in $y$- and in $z$-direction will occur multiple times. This is the only source of inconsistency when using modal truncation for the reduction because the modes themselves do not change.  
\begin{table}[htb]
	\caption{Geometry and material parameters of the 3D cantilever beam and hyperparameter values used for the adaptive sampling.}\label{tab:ParametersBeam}
	\centering
	\begin{tabular}{ c c c c c c } 
		Parameter & Range/Value & Unit & & Hyperparameter & Value \\  \cmidrule{1-3} \cmidrule{5-6}
		Height $h$ & $[0.02, \, 0.05]$ & $\mathrm{m}$ & & $\theta_{lT}$ & $10^\circ$ \\  
		Thickness $t$ & $0.01$ & $\mathrm{m}$ & & $\theta_{uT}$ & $85^\circ$ \\
		Length $l$ & $1.0$ & $\mathrm{m}$ & & $d_{lT}$ & $0.1$ \\
		Young's modulus $E$ & $2.1\cdot 10^{11}$ & $\mathrm{N/m^2}$ & & $d_{uT}$ & $0.2$ \\
		Poisson's ratio $\nu$  & $0.3$ & - & & $d_N$ & 0 \\
		Density $\rho$ & $7860$ & $\mathrm{kg/m^3}$ & & minimum amount of & \multirow{2}{*}{4} \\
		Rayleigh damping $\alpha$ & $8\cdot 10^{-6}$ & $\mathrm{1/s}$ & & samples per cluster  & \\
		Rayleigh damping $\beta$ & $8$ & $\mathrm{s}$ & & &
	\end{tabular}
\end{table}

\begin{figure}[htb]
	\centering
	\includegraphics[width=\textwidth]{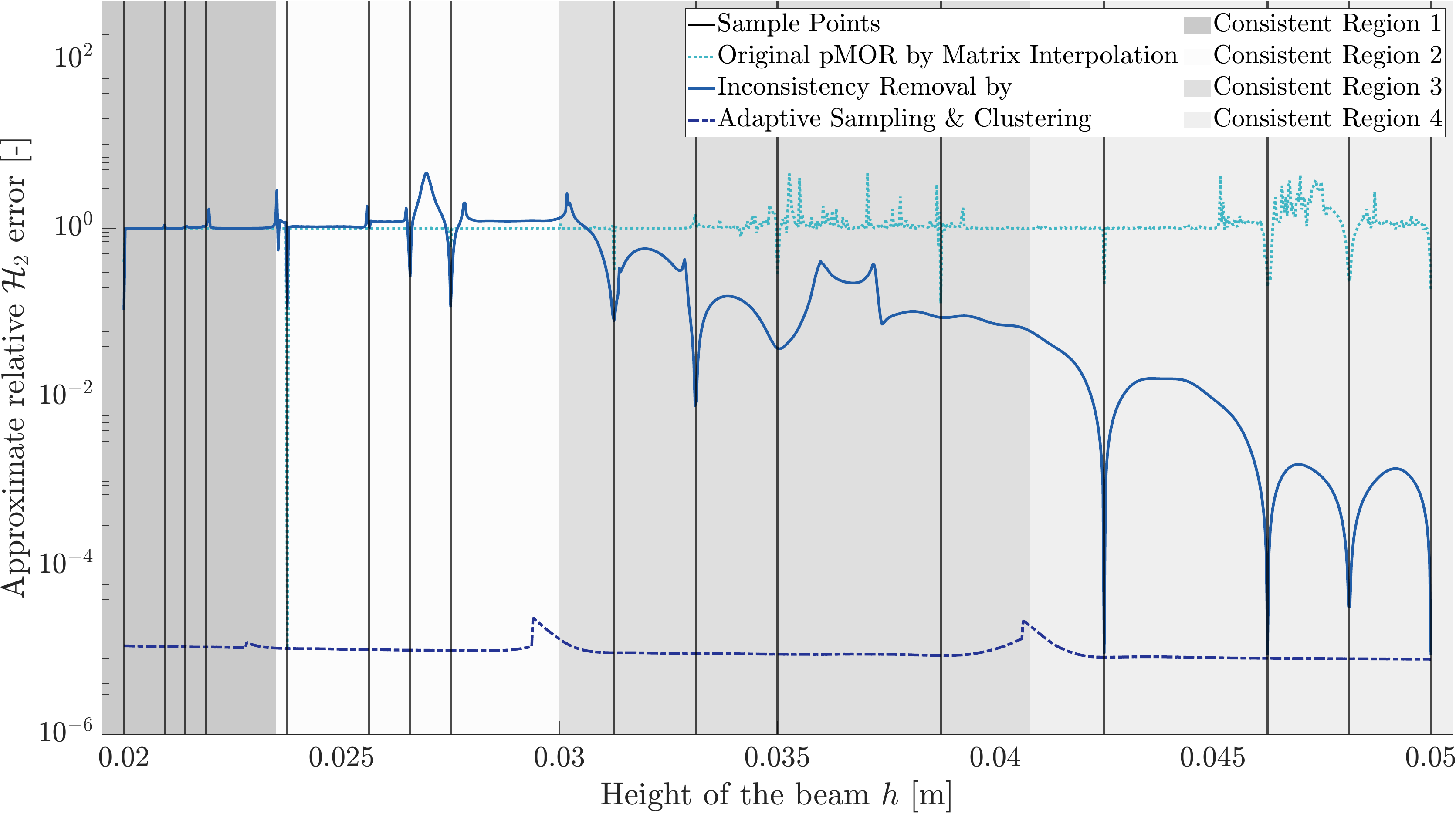}
	\caption{Accuracy of the predicted reduced order models for a cantilever beam.}
	\label{fig:Beam_H2}
\end{figure}

As initial samples, the points $h=2 \, \mathrm{cm}$ and $h = 5 \, \mathrm{cm}$ are used, and the reduced systems are generated by using modal truncation with the first 50 dominant eigenmodes. The hyperparameters of the algorithm are set as shown in \cref{tab:ParametersBeam}. 
The accuracy of the pROM is tested on 1000 linearly spaced sample points in the range from $2 \, \mathrm{cm}$ to $5 \, \mathrm{cm}$. The results of the proposed adaptive sampling and clustering algorithm are compared to the original version of pMOR by matrix interpolation \cite{Panzer2010} and the strategy for removing inconsistencies due to mode switching and truncation by \cite{Amsallem2015}. For all methods, the sample points generated by the adaptive sampling are used. The accuracy of the methods is quantified by the relative $\mathcal{H}_2$ error, which is defined as follows:
\begin{equation}
	\frac{\|H(\cdot; \hat{\p}) - H_r(\cdot ; \hat{\p}) \|_{\mathcal{H}_2}}{\|H(\cdot; \hat{\p})\|_{\mathcal{H}_2}} =  \sqrt{\frac{ \int_{-\infty}^\infty | H(s, \hat{\p}) - H_r(s, \hat{\p}) | \mathrm{d}s }{ \int_{-\infty}^\infty | H(s, \hat{\p}) | \mathrm{d}s }},
\end{equation}
where $H$ and $H_r$ are the transfer functions of the full and the predicted reduced model, respectively. The integral is approximately computed by numerical integration in the above mentioned frequency range. \Cref{fig:Beam_H2} shows the location of the sample points and the approximate relative $\mathcal{H}_2$ error for all test points. \\
The algorithm partitions the parameter space into 4 regions and places 4 samples in each consistent region, which are depicted as differently shaded areas in \cref{fig:Beam_H2}. The interpolated reduced operators that are trained only with information from one consistent region each show a very high accuracy. Especially, in the center of the consistent regions, the interpolation hardly introduces any loss of accuracy. Since the parametric dependency of the full operators on the height of the beam is polynomial with terms up to order three, cubic spline interpolation can capture the dependency exactly. The high accuracy shows that dependencies from the basis have been minimized such that the dependency of the full operators can be learned. Only at the borders between consistent regions, the error increases. This mainly occurs when the classifier does not predict the correct local pROM to use for interpolating the reduced operators, which can, for example, be seen at the border between the consistent regions 2 and 3. Since \textit{k-nearest-neighbor} classification with $k=1$ is used, the border of validity between the two local pROMs is assumed to be in the middle between the last sample point of consistent region 2 and the first sample point of consistent region 3. However, the border between the two consistent regions is actually located at a height of 0.03 m, leading to less accurate results in the area that is misclassified. \\
The original version of pMOR by matrix interpolation suffers from the inconsistency in the training data. Since modal truncation is used as MOR method, inconsistencies due to mode switching and truncation lead to subspace angles of $90^\circ$ between inconsistent bases. Consequently, the product $\mathbf{R}^\intercal \mathbf{V}_k$ becomes singular, which in turn renders the computation of the transformation matrix in \cref{eq:Tk} numerically unstable. Therefore, the transformed reduced operators are singular as well, which explains why the original version of pMOR by matrix interpolation does not even achieve the accuracy of the sampled reduced models at the sample points. \\
The method of inconsistency removal by basis truncation \cite{Amsallem2015} performs slightly better but also shows high errors in most regions. As previously mentioned, mode switching and truncation occurs 3 times leading to 4 consistent regions in this example. The method in \cite{Amsallem2015}, however, only works if mode switching and truncation occurs once. Furthermore, the inconsistency is removed by truncating the basis further which also leads to a loss of accuracy at the sample points.

\subsection{Kelvin Cell \textendash{} 2 Parameters}
\label{sec:KelvinCell_2D}

\begin{figure*}[htb]
	\centering
	\begin{subfigure}[b]{0.47\textwidth}
		\centering
		\includegraphics[width=\textwidth]{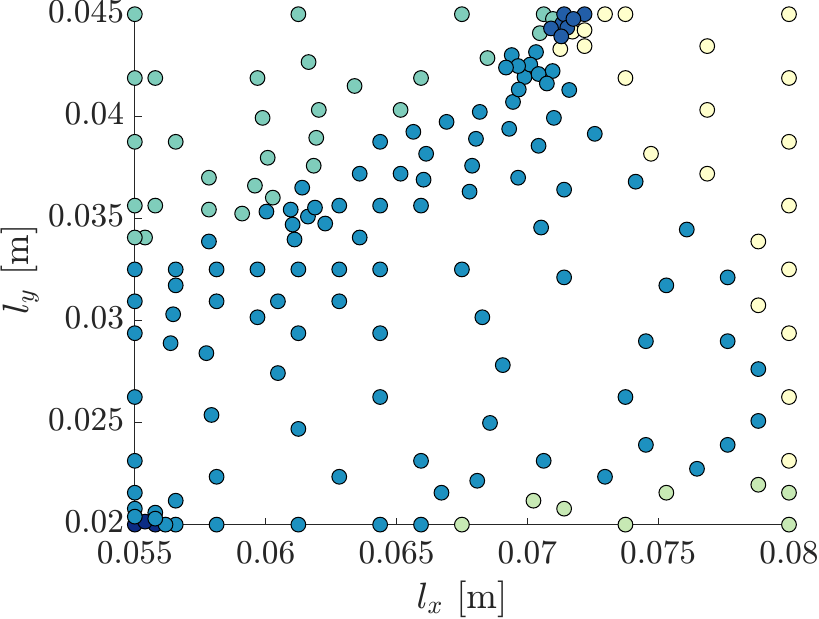}
		\caption{Final distribution and clustering of samples.}  
		\label{fig:SamplesCluster_KelvinCell_2D}
	\end{subfigure}
	\hfill
	\begin{subfigure}[b]{0.51\textwidth}  
		\centering 
		\includegraphics[width=\textwidth]{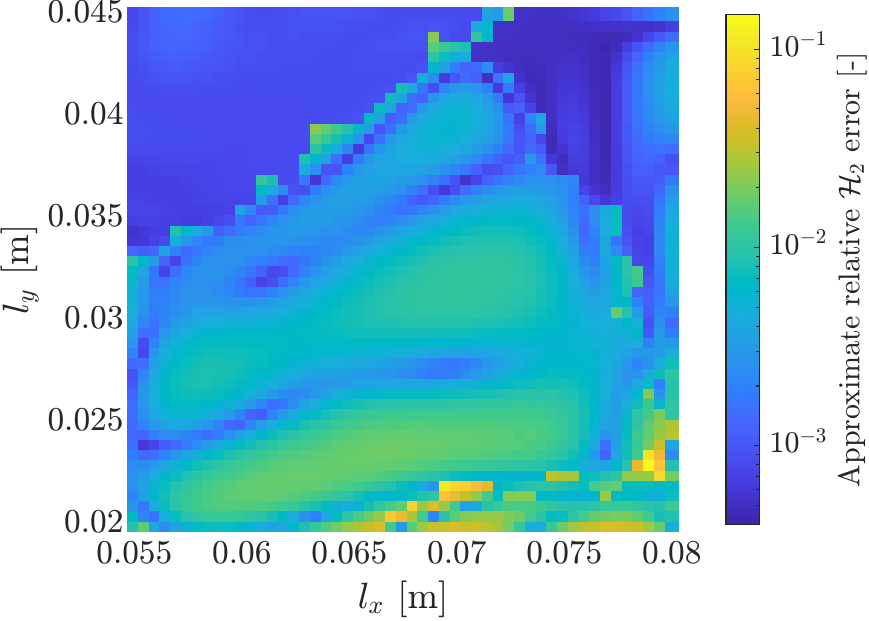}
		\caption{Approximate relative $\mathcal{H}_2$ error of the parametric reduced model.}  
		\label{fig:H2_consistent_KelvinCell_2D}
	\end{subfigure}
	\caption{Results of the adaptive sampling and predicted reduced models for the Kelvin cell with 2 parameters.}  
	\label{fig:H2_KelvinCell_2D}
\end{figure*}

Next, the algorithm is tested on a Kelvin cell structural topology \cite{Thomson*2008} depicted in \cref{fig:KelvinCell}. Twisting, tilting and stretching the cell has been shown to enable controllable anisotropic elasticity of lattice materials \cite{Mao2020}. Each strut is discretized with 50 square Timoshenko beam \cite{Panzer2009} elements. The struts at the bottom are clamped. The input vector $\F$ takes the value 1 at all degrees of freedom (DoFs) of the left face and otherwise 0, the output vector $\G$ takes the value 1 at all DoFs of the right face and otherwise 0. Both, input and output vector, are normalized afterwards. In the first example, the total size of the cell in $x$- and in $y$-direction, $l_x$ and $l_y$, respectively, are used as parameters. The ranges of the parameters and the values for the fixed parameters of the cell are stated in \cref{tab:ParametersKelvinCell}.
\begin{table}[htb]
\caption{Geometry and material parameters of the Kelvin cell and hyperparameter values used for the adaptive sampling.}
\label{tab:ParametersKelvinCell}
	\centering
	\begin{tabular}{ c c c c c c } 
		Parameter & Range/Value & Unit & & Hyperparameter & Value \\  \cmidrule{1-3} \cmidrule{5-6}
		Length $l_x$ & $[0.055, \, 0.080]$ & $\mathrm{m}$ & & $\theta_{lT}$ & $10^\circ$ \\  
		Length $l_y$ & $[0.020, \, 0.045]$ & $\mathrm{m}$ & & $\theta_{uT}$ & $85^\circ$ \\
		Length $l_z$ & $0.05$ & $\mathrm{m}$ & & $d_{lT}$ & $0.1$ \\
		Beam thickness $t$ & $0.001$ & $\mathrm{m}$ & & $d_{uT}$ & $0.2$ \\
		Young's modulus $E$ & $4.35\cdot 10^{9}$ & $\mathrm{N/m^2}$ & & $d_N$ & 0\\
		Poisson's ratio $\nu$  & $0.3$ & - & & $d_C$ & 0.3 \\
		Density $\rho$ & $1180$ & $\mathrm{kg/m^3}$ & & minimum number of & \multirow{2}{*}{1} \\
		Rayleigh damping $\alpha$ & $8\cdot 10^{-6}$ & $\mathrm{1/s}$ & & samples per cluster &  \\
		Rayleigh damping $\beta$ & $8$ & $\mathrm{s}$ & & & 
	\end{tabular}
\end{table}
As initial samples, the 4 corner points of the parameter space are used. The reduced models are generated by modal truncation using the first $50$ eigenmodes as reduced basis. \Cref{tab:ParametersKelvinCell} shows the values used for the hyperparameters of the adaptive sampling. \\
The final distribution and clustering of the samples is shown in \cref{fig:SamplesCluster_KelvinCell_2D}. The algorithm generated 173 sample points and clustered them into 6 regions. These 6 local pROMs are used to predict the transfer function for 1000 linearly spaced frequency points in the interval $[1, \, 1000] \, \mathrm{Hz}$ and for a full grid of $51 \times 51$ test points in the parameter space.

\Cref{fig:H2_consistent_KelvinCell_2D} shows the approximate relative $\mathcal{H}_2$ error for the 2601 test points in the parameter space. In most of the regions, the error of the pROM generated by adaptive sampling and clustering is two to three orders of magnitude lower than for the original version of pMOR by matrix interpolation \cite{Panzer2010} and the strategy for inconsistency removal by \cite{Amsallem2015}, which led to errors of approximately 100 \% in the whole parameter space. However, similar to the results for the cantilever beam, higher errors occur at the borders between inconsistent regions, where the classification model predicts the wrong local pROM to use for predicting the reduced operators. \\
In this example, the offline phase of the adaptive sampling and clustering algorithm took 2969 seconds (wall clock time), where more than 99 \% of that time was required for generating the FE models, assembling the full operators, and computing the reduced operators for the 173 sampled points. This highlights the importance of the adaptive sampling, which aims at placing the sample points optimally with respect to the criteria defined for the adaptive sampling.  Predicting one solution for a queried parameter point in the online phase only takes 0.13 seconds for the pROM. This leads to a tremendous speed-up compared to the FOM, for which generating the FE model, assembling the operators, and solving the full system for a queried parameter point takes 54.4 seconds (wall clock time).

\subsubsection{Detection of inaccuracies at borders}

\begin{figure}[htb]
	\centering
	\label{fig:BorderInconsistentRegions}
	\begin{subfigure}[b]{0.49\textwidth}
		\centering
		\includegraphics[width=\textwidth]{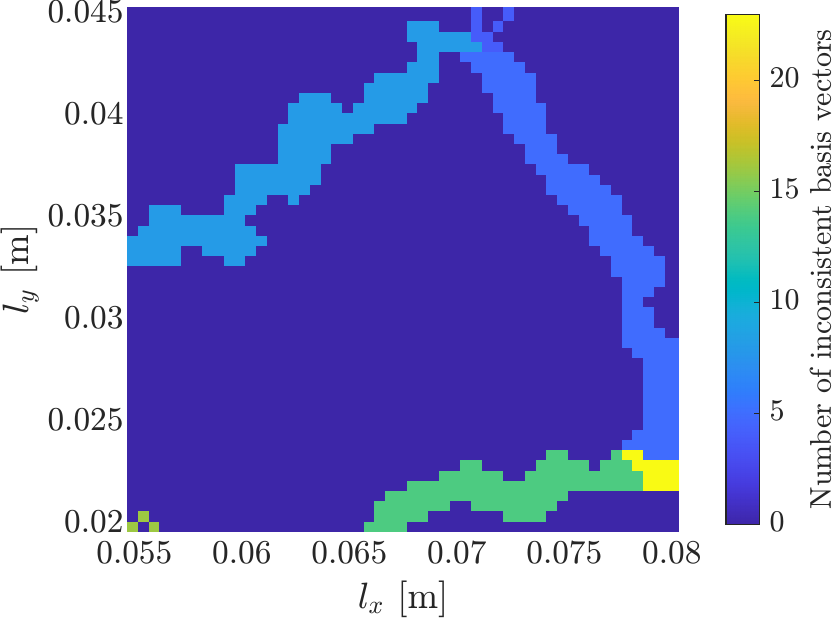}
		\caption{Number of inconsistent basis vectors between the reference bases of the neighboring clusters for each test point.}
		\label{fig:NumberInconsistentBasisVectors}
	\end{subfigure}
	\hfill
	\begin{subfigure}[b]{0.49\textwidth}
		\centering
		\includegraphics[width=\textwidth]{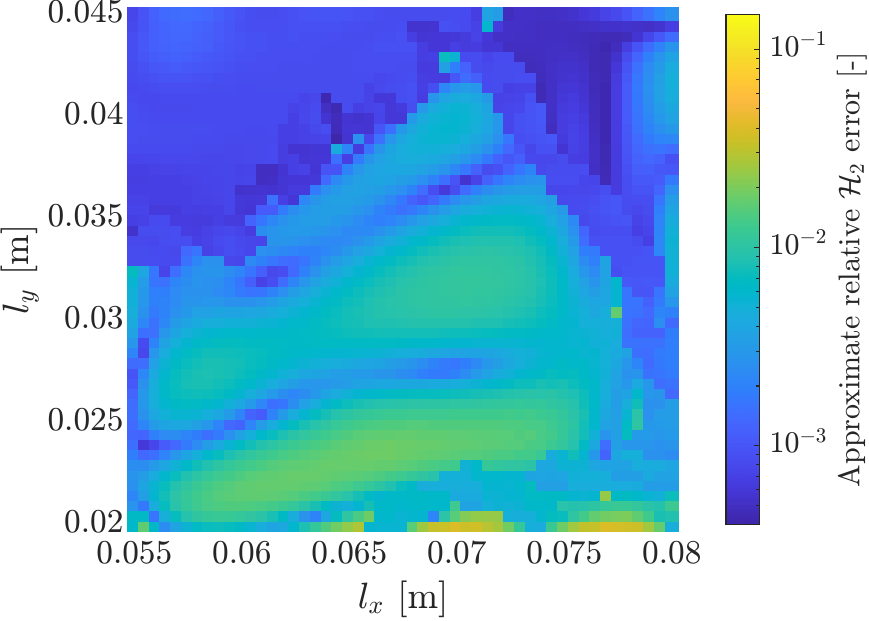}
		\caption{Accuracy of the pROMs by using the proposed global approach for samples between highly inconsistent clusters.}
		\label{fig:H2_RemedyInconsistency}
	\end{subfigure}
 	\caption{Indicator and remedy for inaccurate results at the border between highly inconsistent regions.}
\end{figure}

As shown in the previous examples, the classification model can predict an inappropriate local pROM for a queried parameter point. To evaluate whether the prediction is trustworthy or not, we propose the following indicator for possibly high inaccuracy due to misclassification: We assume that the problem of misclassification only occurs at the border between two inconsistent regions, thus, it only affects queried parameter points that lie between two regions. To identify these points, the triangle of the Delaunay triangulation, in which the queried parameter point lies, is sought. If the samples of this triangle belong to different clusters, it lies between them. In this case, the number of inconsistent vectors can be measured by computing the principal subspace angles between the reference bases of the neighboring inconsistent clusters. The number of subspace angles above the threshold $\theta_{lT}$ can be used to judge whether the prediction is trustworthy or not. This number of subspace angles above the threshold $\theta_{lT}$ is shown in \Cref{fig:NumberInconsistentBasisVectors} for each test point from the comparison in \cref{sec:KelvinCell_2D}. A comparison to the approximate relative $\mathcal{H}_2$ error depicted in \cref{fig:H2_consistent_KelvinCell_2D} shows that the borders between regions are correctly identified and the number of inconsistent basis vectors correlates well with the error in case of misclassification. Note that this is merely an indicator for areas of decreased accuracy and cannot be used as an error estimate. However, it avoids assembling the full order model in the complete parameter space, which would be required for an error estimate, and is thus computationally substantially cheaper.\\ 
To achieve better results in these regions, a global pMOR approach can be used: For each queried parameter point, the reduced bases of the neighboring samples identified by the triangulation can be concatenated to a global basis. Next, the FOM is assembled for this parameter point and projected onto this global basis. Generating this ROM takes more time than with pMOR by matrix interpolation since the FOM must be assembled, but it results in a high accuracy. This is shown in \cref{fig:H2_RemedyInconsistency}, where this global pMOR approach is used for all test points that lie between inconsistent regions, and the approach of pMOR by matrix interpolation is used for all other test points.

\subsection{Kelvin Cell \textendash{} 3 parameters}

\begin{figure}[htb]
	\centering
	\begin{subfigure}[b]{0.49\textwidth}
		\centering
		\includegraphics[width=\textwidth]{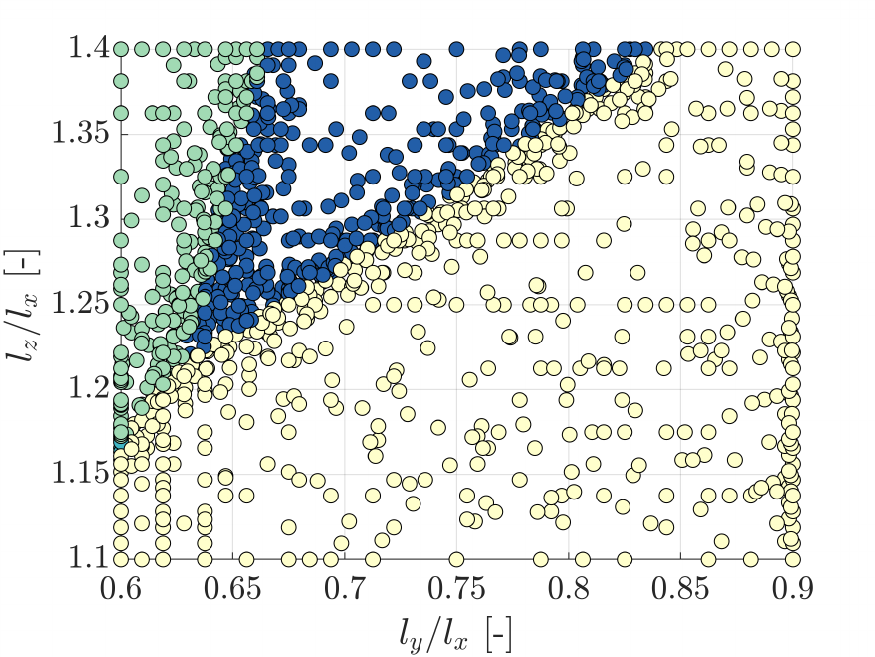}
		\caption{Final distribution and clustering of samples.}
		\label{fig:KelvinCell3D_Samples}
	\end{subfigure}
	\hfill
	\begin{subfigure}[b]{0.49\textwidth}
		\centering
		\includegraphics[width=\textwidth]{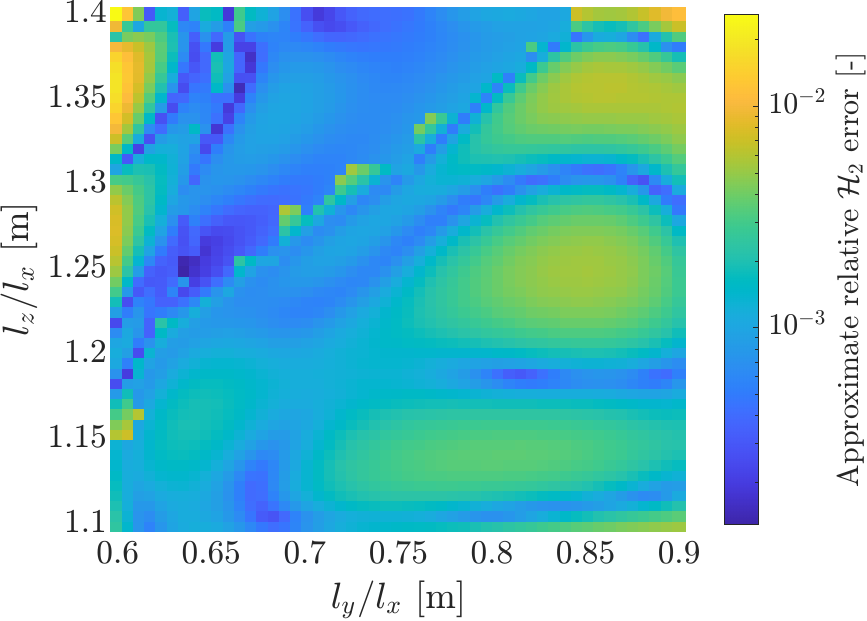}
		\caption{Accuracy of the predicted reduced models for $l_x = 50$ mm.}
		\label{fig:KelvinCell3D_H2}
	\end{subfigure}
	\caption{Results of the adaptive sampling and clustering for the Kelvin cell with 3 parameters.}
	\label{fig:KelvinCell3D}
\end{figure}

For an example with a higher-dimensional parameter space, the Kelvin cell structural topology (\cref{fig:KelvinCell}) is investigated again. Now, $l_x$ and the ratios $l_y/l_x$ and $l_z/l_x$ are used as parameters and varied within the intervals $[45, \, 55]$ mm for $l_x$, $[0.6, \, 0.9]$ for $l_y/l_x$ and $[1.1, \, 1.4]$ for $l_z/l_x$. The fixed parameters of the model and the hyperparameters of the algorithm stay mostly the same as denoted in \cref{tab:ParametersKelvinCell}. Only $d_N$ is now set to 0.05 in order to prevent the algorithm from placing many new samples in small regions where skew edges occur. As initial samples, again only the corner points are used. MOR is performed by modal truncation using the first 50 eigenmodes. \\
The adaptive sampling results in 1389 samples that are clustered into 4 groups. \Cref{fig:KelvinCell3D_Samples} shows the distribution of the samples projected in the $l_y/l_x$-$l_z/l_x$-plane. The accuracy of the predicted pROMs is assessed for 1000 linearly spaced frequency points in the interval $[1, \, 1000]$ Hz and a full grid of $21 \times 21 \times 21$ test points in the parameter space. The mean approximate relative $\mathcal{H}_2$ error for these test points results in $0.21 \%$, whereas the original version of pMOR by matrix interpolation \cite{Panzer2010} and the strategy for inconsistency removal by \cite{Amsallem2015} again led to errors of approximately 100 \% in the whole parameter space. For a more detailed analysis of the accuracy of the pROM, the approximate relative $\mathcal{H}_2$ error is shown in the $l_y/l_x$-$l_z/l_x$-plane for $l_x = 50$ mm in \cref{fig:KelvinCell3D_H2}. As for the example with two parameters, the error increases between the borders of the clusters and in less densely sampled regions. 

\section{Discussion \& Conclusion}
\label{sec:conclusions}

\subsection{Discussion}

The proposed algorithm has some limitations, which have already been discussed partially. First of all, the accuracy at borders between highly inconsistent regions might deteriorate if the queried point is misclassified. However, by using a global pMOR approach at borders between inconsistent regions, a better accuracy can be achieved at the cost of a higher computational effort.  \\
Secondly, pMOR by matrix interpolation in general, suffers from the curse of dimensionality. This is because, in higher dimensions, many sample points are required to obtain a good interpolation. Additionally, the adaptive sampling becomes computationally more complex with increasing dimension of the parameter space. However, compared to an adaptive sampling on a full grid, the proposed framework is computationally more efficient as fewer samples are required in each iteration of the adaptive sampling because no regular grid must be maintained. \\
For applications yet to be identified, the dynamics of the system might change so strongly in the whole parameter space that also for a fine sampling, all samples have a different local reduced basis, and a clustering into larger consistent regions is not possible. In this case, other pMOR methods should be used.

\subsection{Conclusion}

A general framework for removing inconsistencies in the training data of pMOR by matrix interpolation is presented. These inconsistencies stem from strong changes in the sampled reduced bases over the parameters, which can occur due to the MOR method used, a change of the system dynamics with change of the parameters, and mode switching and truncation. By using modal truncation for the reduction, performing an adaptive sampling and partitioning the parameter space into several consistent regions, all of these sources of inconsistencies are minimized within a unified framework. The proposed algorithm is applied to a cantilever Timoshenko beam and the Kelvin cell topological structure for one-, two- and three-dimensional parameter spaces. Comparisons to the original version of pMOR by matrix interpolation and a different strategy to remove inconsistencies due to mode switching and truncation show that the proposed algorithm leads to errors that are several orders of magnitude lower at only slightly higher computational effort for low-dimensional parameter spaces.

\addcontentsline{toc}{section}{References}
\bibliographystyle{plainurl}
\bibliography{bibtex/references}

\end{document}